\documentclass[11pt]{article}
\catcode`\@=11 \catcode`\@=12

\usepackage{amsmath}
\usepackage{amssymb}
\usepackage{amsthm}
\usepackage{oldgerm}

\newcommand{\Rm}{\mathbb{R}}

\newcommand{\Ru}{\ensuremath{\mathfrak{R}}}
\newcommand{\Cu}{\ensuremath{\mathfrak{C}}}

\newcommand{\mN}{\ensuremath{\mathcal{N}}}

\newcommand{\Nm}{\ensuremath{\mathbb{N}}}
\newcommand{\Zm}{\ensuremath{\mathbb{Z}}}
\newcommand{\Lm}{\ensuremath{\mathbb{L}}}

\newcommand{\mK}{\ensuremath{\mathcal{K}}}

\newcommand{\mA}{\ensuremath{\mathcal{A}}}
\newcommand{\mO}{\ensuremath{\mathcal{O}}}

\newcommand{\mI}{\ensuremath{\mathcal{I}}}

\newcommand{\Tm}{\ensuremath{\mathbb{T}}}
\newcommand{\G}{\ensuremath{\overline {\mathcal{G}}}}

\newcommand{\vs}{\vspace{.2cm}}

\newtheorem{lem}{Lemma}
\newtheorem{thm}{Theorem}

\newtheorem{cor}[lem]{Corollary}
\newtheorem{prop}[lem]{Proposition}
\newtheorem{defn}[lem]{Definition}

\def\proof {\noindent{\sc{Proof. }}}
\def\qed {\mbox{}\hfill {\small \fbox{}} \\}
\def\lto{\longrightarrow}
\def\lmto{\longmapsto}
\def\eq{\Longleftrightarrow}
\def\leq{\leqslant}
\def\geq{\geqslant}

\title{On the Conley decomposition of Mather sets}
\author{Patrick  Bernard}
\begin{document}

\maketitle
\begin{center}
-----
\end{center}

\begin{small}
\vs \noindent
 Patrick Bernard
 \footnote{membre de l'IUF}\\
Universit\'e  Paris-Dauphine\\
CEREMADE, UMR CNRS 7534\\
Pl. du Mar\'echal de Lattre de Tassigny\\
75775 Paris Cedex 16,
France\\
\texttt{patrick.bernard@ceremade.dauphine.fr}\\
\end{small}
\vs
\thispagestyle{empty}
\begin{center}
-----
\end{center}
\vs
\noindent
\textbf{R\'esum\'e :}
Dans le contexte de la th\'eorie de Mather des syst\`emes
Lagrangiens, on \'etudie la d\'ecomposition en composantes 
transitives par chaines des ensembles invariants de Mather.
Comme application, on montre, sous des hypoth\`eses appropri\'ees,
la semi-continuit\'e de l'ensemble d'Aubry.
\vs

\begin{center}
-----
\end{center}

\vs\noindent
\textbf{Abstract : }
In the context of Mather's theory of Lagrangian systems,
we study the decomposition in chain-transitive classes of the 
Mather invariant sets. As an application, we prove,
under appropriate hypotheses, the semi-continuity
of the so-called Aubry set as a function of the Lagrangian.

\begin{center}
-----
\end{center}

\vspace{1cm}


\newpage

In the study of Lagrangian systems, John Mather introduced 
several invariant sets composed of globally minimizing extremals.
He developed methods to construct several orbits undergoing
interesting behaviors in phase space under some assumptions 
on these invariant  sets, see \cite{Ma:93}.
In order to pursue this theory and to apply it on examples, it is
necessary to have tools to describe precisely the 
invariant sets.
At least two points of view can be adopted.
One can study the invariant set from a purely topological
point of view 
in the style of Conley
as compact metric spaces with flows, and study
their transitive components.
One can also study these set from the point of view of action 
minimization, and  decompose them in invariant subsets
that have been called static classes.
These points of view are very closely related, but each of 
them has specific features. For example, understanding the 
decomposition in static classes is necessary for the variational
construction of interesting orbits, while the topological 
decomposition behaves well under perturbations.

Our goal in the present paper is to explicit the links
between these two decompositions. We explain that the topological
decomposition is finer than the variational one, and that they
coincide for most (but not all) systems.
As an application, we prove a result of semi-continuity
of the so-called Aubry set as a function of the Lagrangian, 
under certain non-degeneracy hypotheses.
The semi-continuity of the Aubry set is a subtle problem, 
which has remained open for several years, until
John Mather gave a counter example,
see \S 18 in \cite{Ma:04}.
In the same paper, he also states without proof 
that semi-continuity holds under appropriate hypotheses.
Our result extends the one of Mather.
The methods we use are inspired from the recent work of Fathi, Figalli
and Rifford, \cite{FFR}.

\section{Introduction}
We have to recall the fundamental constructions of Mather
theory before we can state our results.
We rapidly expose them without proofs. These proofs are available
in 
\cite{Ma:91,Ma:93,Fa:book,Fourier,AMS}.

We consider a compact manifold without boundary $M$.
It is convenient to endow once and for all this manifold with
a Riemannian metric. 
We shall work in the standard framework of Mather theory,
and study  $C^2$ Lagrangians 
$L(t,q,v):\Tm\times TM\lto \Rm$, where $\Tm=\Rm/\Zm$.
Given $t\in \Rm$, we will also denote by $t$ the associated 
element of $\Tm$. Given $\tau\in \Tm$, we will also denote
by $\tau$ the corresponding element of $[0,1[$.
The main object of study is the dynamics of minimizing extremals 
of $L$.
An absolutely continuous curve $q(t):I\lto M$ is called a minimizing extremal if,
for each $t_0\leq t_1$ in $I$ and each absolutely continuous curve
$\gamma(t):[t_0,t_1]\lto M$,  satisfying $\gamma(t_0)=q(t_0)$
and $\gamma(t_1)=q(t_1)$, we have
$$
\int_{t_0}^{t_1}L(t,q(t),\dot q(t))dt \leq
\int_{t_0}^{t_1}L(t,\gamma(t),\dot \gamma(t))dt.
$$
The curve $q(t)$ is called a locally minimizing extremal if
each time $t$ is contained in the interior of an interval $J$ such
that $q$ is a minimizing extremal on $J$.
We assume \vs

\noindent
\textbf{Convexity :} For each $(t,x)\in \Tm\times M$, the function
$v\lto L(t,x,v)$ is convex, has positive definite Hessian at each point, and is superlinear. In short, for each $(t,x)\in \Tm\times M$,
we have
$$
\partial^2_{v}L(t,x,v)>0,\forall v\in T_xM\quad 
\text{and}\quad
\lim_{\lambda \rightarrow \infty}\big(L(t,x,\lambda v)/\lambda\big)
=\infty \quad\forall v\in T_xM-\{0\}.
$$ 

Under the convexity hypothesis,  there exists a vector-field
$E_L$ on $\Tm\times TM$, the Euler-Lagrange vector-field,
 such that a curve $q(t)$ 
is a local minimizing extremal if and only if the associated curve
$t\lmto (t, q(t), \dot q(t))$ is an integral curve of $E_L$.
We assume \vs

\noindent
\textbf{Completeness :} The Euler-Lagrange vector-field 
has a complete flow $\phi^t$ on 
$\Tm\times TM$.\vs

A $C^2$ Lagrangian satisfying convexity and completeness will 
be called a Tonelli Lagrangian in the sequel.
It is useful to define
the function $A_L:\Rm\times M\times \Rm\times M\lto \Rm$ by the expression
$$
 A_L(S,q;T,r):= \min_{\gamma(S)=q, \gamma(T)=r} 
\int_S^TL(t,\gamma(t),\dot \gamma(t))dt,
$$
where the minimum is taken on the set of $C^1$ curves
$\gamma(t):[S,T]\lto M$ 
which satisfy $\gamma(S)=q$ and $\gamma(T)=r$.
The existence of the minimum for a Tonelli Lagrangian
is a standard result derived from  Tonelli's work.
It is known that there exists a unique constant
$\alpha(L)$ such that the function 
$
\tilde A_L(S,q,T,r)+(T-S)\alpha(L)
$
is bounded on $\{T\geq S+1\}$.
This constant  is sometimes called the Ma\~n\'e critical value,
although it was first introduced by Mather in $\cite{Ma:91}$.
Most of the dynamics of locally minimizing orbits
is encoded in the function $h_L:\Tm\times M\times \Tm\times M\lto \Rm$
defined, following Mather, by
$$
h_L(\tau,q;\theta,r):=\liminf_{\Nm \ni T \lto \infty} 
\big( A_L(\tau,q;\theta+T,r)+ (T+\theta-\tau)\alpha(L)\big).
$$
The function $h_L$  is Lipschitz
continuous (and semi-concave) on $\Tm\times M\times \Tm \times M$.
A function $u:\Tm\times M\lto \Rm$ is said dominated by $L$ if 
\begin{equation}\tag{D}
u(t, q(t))-u(s, q(s))\leq \int_s^t L(\sigma,q(\sigma),\dot q(\sigma))
+\alpha(L) d\sigma
\end{equation}
for each curve $q(\sigma)\in C^1(\Rm, M)$ and each $s<t$ in $\Rm$.
This implies that 
$$
u(y)-u(x)\leq   
h_L(x,y)\qquad \forall x,y\in\Tm\times M.
$$

As was noticed by Albert Fathi, the relevance of dominated functions
is that there are invariant sets of the Euler-Lagrange flow naturally
associated to them.
In order to define these sets, it is necessary first to define,
following Fathi,
the notion of calibrated curve.
A curve $q(t): I\lto M$ is  said calibrated by the dominated 
function $u$ if, for each $s<t$ in $I$, the equality holds in (D).
It is clear that calibrated curves are minimizing extremals.
For each dominated function $u$, we define the set
$$
\tilde \mI(L,u)\subset \Tm\times TM
$$
as follows: 
$\tilde \mI(L,u)$ is the set of points $(\tau,q,v)$ such that 
there exists a calibrated curve $q(s):\Rm\lto M$ 
satisfying  $(\tau,q,v)=(\tau, q(\tau),\dot q(\tau))$.
It is known that 
$\tilde \mI(L,u)$ is a compact invariant set
of the Euler-Lagrange flow.
The projection $\pi:(t,x,v)\lmto(t,x)$ induces a bi-Lipschitz
homeomorphism between $\tilde \mI(L,u)$ and its image
$\mI(L,u)\subset \Tm\times M$.
We shall always endow $\tilde \mI(L,u)$ with the flow
induced from the Euler-Lagrange flow, and $\mI(L,u)$
with the conjugated flow.

The function $u$ is called a Weak KAM solution if 
it is dominated and if, in addition,
for each point $(\tau,q)\in \Tm\times M$, there exists
a   calibrated curve 
$q(s):(-\infty, \tau]\lto M$ such that  $q(\tau)=q$.
Given a Weak KAM solution $u$, we define the set 
$$
\G(L,u)\subset \Tm\times TM
$$
as the set of points $(\tau,q,v)$ such that there exists 
 a   calibrated curve 
$q(s):(-\infty, \tau]\lto M$ satisfying $q(\tau)=q$
and $\dot q(\tau)=v$.
The set $\G(L,u)$ is compact and negatively invariant, so that
the Euler-Lagrange flow defines on it a  semi-flow
(with negative times). Note also that $\pi(\G(L,u))=\Tm\times M$ and
 that 
$$
\tilde \mI(L,u)=\bigcap _{t\leq 0} \phi^t\big(\G(L,u)\big).
$$
It is known that, for each $x_0\in M$, the function 
$x\lmto h_L(x_0,x)$ is a Weak KAM solution.

The set
$$
\tilde \mA(L):= \bigcap \tilde \mI(L,u)
$$
is not empty, where the intersection is taken on the set 
of all dominated functions, or equivalently on the set of all 
Weak KAM solutions. This is the definition of Fathi
of a set previously introduced by John Mather
in \cite{Ma:03}, and called the Aubry 
set.
It is clearly compact and invariant, we shall always endow
it with the Euler-Lagrange flow. The projection $\pi$
restricted to $\tilde \mA(L)$ is a bi-Lipschitz homeomorphism
into its image $\mA(L)$. We endow this image with the conjugated flow.
Let us recall from now on that Fathi proved the existence of a Weak KAM
solution $u$ such that $\mI(L,u)=\mA(L)$.
John Mather noticed that the 
function 
$$d_L(x,y):=h_L(x,y)+h_L(y,x)
$$
is a pseudo-metric on $\mA(L)$. 
Indeed, it is known that $d_L(x,y)\geq 0$ for each $x$ and $y$
in $\Tm\times M$, and that 
$d_L(x,x)=0$ if and only if $x$ belongs to $\mA(L)$.
The function $d_L$ 
 is symmetric and satisfies the triangle inequality,
but there may exist points $x\neq y$ such that $d_L(x,y)=0$.
The relation $\Ru$ on $\mA$ defined by
$$
x \Ru y \Longleftrightarrow
d_L(x,y)=0
$$
 is then an equivalence relation
on $\mA(L)$. The equivalence classes are called the static classes.
They are compact invariant subsets of $\mA(L)$.
Note that the pseudo-metric $d_L$ descends to a metric,
that we still denote by $d_L$, on the set $\dot \mA(L)$ of
static classes. The set $\dot\mA(L)$, endowed with the metric 
$d_L$ has been called by Mather the quotient Aubry set.
It is a compact metric space.

The set 
$$
\tilde \mN(L):= \bigcup \tilde \mI(L,u)
$$
is  compact and invariant, where the union  is taken on the set 
of all dominated functions, or equivalently on the set of all 
Weak KAM solutions. This is the definition of Fathi
of a set previously introduced by John Mather
\cite{Ma:03}, and called the Ma\~n\'e
set.
The $\alpha$ and $\omega$-limits of orbits of the Ma\~n\'e set
are contained in the Aubry set, see for example \cite{Fourier}
or \cite{Ma:03}.
More is true: given an orbit of the Ma\~n\'e set, there exists
a static class which contain all its $\alpha$-limit points
and a static class which contain all its
$\omega$-limit points. These static classes are equal if and 
only if the orbit is contained in the Aubry set.

For each Weak KAM solution $u$, we define  the relation
$x \Ru_u y$  on $\Tm\times M$ by
$$
x \Ru_u y \Longleftrightarrow u(y)-u(x)=h_L(x,y).
$$
This relation is  transitive.
Indeed, we have 
$
x \Ru_u y \Longleftrightarrow u(y)-u(x)\geq h_L(x,y)
$
(the converse inequality always holds).
If $x_0\Ru_u x_1$ and $x_1 \Ru_u x_2$,
then we have
$$
u(x_2)-u(x_1) \geq h_L(x_0,x_1)+h_L(x_1,x_2)\geq h_L(x_0,x_2)
$$
so that $x_0 \Ru_u x_2$.
The Aubry set $\mA(L)$ is the set 
of points $x\in \Tm\times M$ such that $x\Ru_u x$.
The symmetrized relation is nothing but $\Ru$: 
$$
\big(x\Ru_u y \;\text{and}\;
y\Ru_u x\big)
\Longleftrightarrow
\big(x\in \mA(L)\;\text{and}\;y\in \mA(L)\;\text{and}\; d_L(x,y)=0\big).
$$

We  denote by $\Cu_u$ the  relation of chain-connection  
on $\mI(L,u)$, see Appendix. We use the same symbol for 
the relation of chain connection in
 $\tilde \mI(L,u)$ and in  $\G(L,u)$
 (these relations coincide on $\tilde \mI(L,u)$ by Lemma \ref{add}
 of the Appendix).

\begin{prop}\label{oneway}
For each Weak KAM solution $u$, and for $x$ and $y$ in $\mI(L,u)$,
we have
$$
x\Ru_u y \Longrightarrow x\Cu_u y.
$$
More precisely, given two points $x$ and $y$ in $\Tm\times M$,
the relation $x\Ru_u y$ implies that $x\in \mI(L,u)$ and that there
exists points $\tilde x$ and $\tilde y$ above $x$ and $y$
in $\G(L,u)$ such that $\tilde x\Cu_u \tilde y$. 
\end{prop}

The proof will be given in Section \ref{convergence}.
The following statement is due tu Ma\~n\'e, see \cite{Mn:97},
 and a proof can be found
in \cite{Fa:book}. 

\begin{cor}
The Ma\~n\'e set $\tilde \mN(L)$ is chain transitive.
The Aubry set $\mA(L)$ is chain recurrent.
Each static class is chain transitive in $\mA(L)$.
\end{cor}

We give a proof as an  application of Proposition \ref{oneway}.

\proof
In order to prove the chain recurrence of the Aubry set,
let us first recall that there exists a 
weak KAM solution $u$ 
such that $\mI(L,u)=\mA(L)$.
For this function $u$, we have 
$x\Ru_u x\Rightarrow x\Cu_u x$, which precisely
says that the point $x$ is chain-recurrent in $\mA(L)$.
If $x$ and $y$ belong to the same static
class, then we have, still with the same function $u$, that
 $x\Ru_u y$, and therefore 
$x\Cu_u y$.
Finally, the chain transitivity of the Ma\~n\'e set can be proved
in several steps.
First, let $x$ and $y$ be points of $\mA(L)$.
Then, taking $u=h(x,.)$, we see that $x\Ru_u y$,
and therefore $x\Cu_u y$. As a consequence, $x$ and $y$ can 
be connected by chains in $\mI(L,u)$, 
and therefore the points $\tilde x$ and $\tilde y$ which
are the points of $\tilde \mI(L,u)$ above $x$ and $y$ 
can be connected by chains in $\tilde \mN(L)$.
Now since every point $\tilde x\in \tilde \mN(L)$
has its $\omega$-limit contained in $\tilde \mA(L)$, 
there exists a point $\tilde \omega \in \tilde \mA(L)$
such that $\tilde x\Cu_{\tilde \mN} \tilde\omega$.
In the same way, for each $\tilde y\in \tilde \mN(L)$, there exists
a point $\tilde \alpha \in \tilde \mA$ such that 
$\tilde \alpha \Cu_{\tilde \mN} y$.
By transitivity of the relation 
$\Cu_{\tilde \mN}$, we conclude that $\tilde x\Cu_{\tilde \mN}\tilde y$.
\qed

It is known that the converse to Proposition \ref{oneway}
does not hold in general. However, it holds in many examples,
and has interesting consequences, that we now describe.

\begin{defn}
We say that the Lagrangian $L$ satisfies the coincidence hypothesis
if, for each  Weak KAM solution $u$, the relations $\Cu_u$ and 
$\Ru_u$ coincide on $\mI(L,u)$.
\end{defn}

By well-known properties of the Conley decomposition
of flows recalled in the Appendix, we obtain:

\begin{prop}
If $L$ satisfies the coincidence hypothesis,
 then, for each Weak KAM solution $u$,
the Aubry set is precisely the chain-recurrent set of $\G(L,u)$.
The static classes  are the chain transitive components
 of $\mA(L)$, they are also the connected components of $\mA(L)$.
  The quotient Aubry set is totally disconnected. 
\end{prop}

The coincidence hypothesis also has as a consequence
the semi-continuity of the Aubry set.
In order to be more precise, we now introduce a notion
of convergence for Tonelli Lagrangians.
The sequence $L_k$ of Tonelli Lagrangians is said to converge to $L$
if $\{L_k\}$ is a uniform family of Tonelli Lagrangians
(see definition in Section \ref{convergence}),
and if $L_k$ converge to $L$ uniformly on compact sets
as $k\lto \infty$. We shall study this convergence in 
Section \ref{convergence}.

\begin{thm}\label{semi}
Let $L$ be a Tonelli Lagrangian satisfying the coincidence hypothesis.
Let $L_k$ be a sequence of Tonelli Lagrangians converging 
to $L$. Let $U$ be a neighborhood of $\tilde \mA(L)$ in 
$ \Tm\times TM$.
Then, there exists $k_0$ such that $\tilde \mA(L_k)\subset U$
for each $k\geq k_0$.

In general (without the coincidence hypothesis),
the set $\limsup \tilde \mA(L_k)$ is contained in the chain recurrent
set of $\tilde \mN(L)$ (and more precisely in the union 
of the chain recurrent sets of $\tilde \mI(L,u)$). 
\end{thm}

John Mather has stated this result without proof  in
\cite{Ma:04} under the hypothesis that the quotient Aubry set 
$\dot \mA(L)$
has vanishing $1$-dimensional Hausdorff measure. 
We shall see that this hypothesis
of Mather implies the coincidence 
hypothesis hence the statement above is stronger than Mather's one.
A very partial version of Theorem \ref{semi}
 was proved earlier in \cite{ETDS}.
The examples described in \cite{Ma:04} show that the semi-continuity
of the Aubry set with respect to the Lagrangian is not always true.

In order to extract more information from
the constructions presented above, Mather noticed that,
when $\omega$ is a closed one-form on $ M$,
the Lagrangian
$$(L+\omega)(t,x,v):=L(t,x,v)+\omega_{x}(v)$$
has the same Euler-Lagrange extremals than $L$, but not
the same Aubry set.
The Aubry set of $L+\omega$ depends only on the De Rham 
cohomology class
of $\omega$ in $H^1(M,\Rm)$. 
Given $c\in H^1(M,\Rm)$, we denote  by $\mA(L+c)$
the Aubry set of $L+\omega$, where $\omega$ is any closed
form of cohomology $c$.
The semi-continuity of the Aubry set
 with respect to the cohomology $c$
is an open question, but we have:

\begin{cor}
Given a Tonelli Lagrangian $L$, the set-valued map
$$H^1(M,\Rm)\ni c\lmto \tilde \mA(L+c)
$$
is upper-semicontinuous at every point $c$ where $L+c$ satisfies
the coincidence hypothesis.
If $L$ is generic in the sense of Ma\~n\'e, then this hypothesis
holds for each $c$.
\end{cor}

In the sequel, we detail and prove what has been stated above.
All is based on basic properties of  the Conley decomposition 
of flows which are recalled in the Appendix.
In Section \ref{convergence}, we study the main features of 
the convergence of Tonelli Lagrangians.
We prove Proposition \ref{oneway} and Theorem \ref{semi}.
Finally, in Section \ref{coincidence}, we discuss the hypothesis of coincidence.

\textbf{Acknowledgements :}
I thank Albert Fathi for pointing out  that
$5\Rightarrow 4$ in Section \ref{coincidence}.

\section{Convergence of Tonelli Lagrangians}\label{convergence}
We define the notion of convergence that is used throughout
the paper, and collect its relevant properties.
It is first useful to recall from \cite{AMS} the notion
of uniform family of Tonelli Lagrangians.
A Family $\Lm\subset C^2(\Tm\times TM,\Rm)$ of Tonelli Lagrangians 
is called uniform if:
\begin{itemize}
\item[$(i)$]
There exist two superlinear functions $l_0$ and $l_1:\Rm^+\lto \Rm$
such that each Lagrangian  $L$ of the family satisfies 
$l_0(|v|_q)\leq L(t,q,v)\leq l_1(|v|_q).$
\item[$(ii)$]
There exists an increasing function $K(k): \Rm^+\lto \Rm^+$ such 
that, if   $\phi $ 
is the Euler-Lagrange flow of a Lagrangian of the Family,
then, for each $t\in [-1,1]$, we have 
$$
\phi^t\big( \{|v|_q\leq k\}\big)
\subset \{|v|_q\leq K(k)\}\subset \Tm \times TM.
$$
\item[$(iii)$]
There exists a finite atlas $\Psi$ of 
$M$ such that, for each chart $\psi\in \Psi$ and each
Lagrangian $L$ of the family,
we have $\|d^2(L\circ T\psi)_{(t,q,v)}\|\leq K(k)$
for $|v|_q\leq k$.
\end{itemize}
The following result  is proved in the Appendix of \cite{AMS}:

\begin{prop}
The set of functions $u$ which are Weak KAM solutions of some $L\in\Lm$
is equi-Lipschitz. 
For each compact interval $[a,b]$, the set of   curves 
$q(s):[a,b]\lto M$
which are minimizing extremals for some $L\in \Lm$ 
is relatively compact in $C^1([a,b],M)$.
\end{prop}

The sequence $L_k$ of Tonelli Lagrangians is said to converge to $L$
if it forms a uniform family and if $L_k$ converges to $L$ uniformly
on compact sets.

\begin{lem}
If $L_k\lto L$, then $\alpha(L_k)\lto \alpha(L)$.
If $u_k$ is a sequence of Weak KAM solutions of $L_k$
which converge uniformly to $u$,
then $u$ is a Weak KAM solution of $L$.
\end{lem}
\proof
Let $u_k$ be a sequence of Weak KAM solution of $L_k$.
Let us take a subsequence such that $\alpha(L_k)$ has a limit
$\alpha$ and such that $u_k$ has a uniform limit $u$.
Taking the limit in the inequality
$$
u_k(t,q(t))-u_k(s,q(s))\leq \int_s^t
 L_k(\sigma, q(\sigma),\dot q(\sigma)) +\alpha_k(L) d\sigma
$$
we get 
$$
u(t,q(t))-u(s,q(s))\leq \int_s^t
L(\sigma, q(\sigma),\dot q(\sigma)) +\alpha d\sigma
$$
for each curve $q(\sigma)\in C^1([s,t],M)$.
This implies that $\alpha\geq \alpha(L)$.
For each $x=(\tau,q)\in \Tm\times M$, let
 us now consider a sequence $q_k\in C^1((-\infty,\tau],M)$
such that $q_k(\tau)=q$ and such that
$$
u_k(\tau,q_k(\tau))-u_k(s,q_k(s))=\int_s^{\tau}
 L_k(\sigma, q_k(\sigma),\dot q_k(\sigma)) +\alpha_k(L) d\sigma
$$
for each $s<\tau$.
This sequence has a subsequence which converges in 
$C^1_{loc}((-\infty,\tau],M)$
to a limit $q(t)$, which satisfies
$$
u(\tau,q(\tau))-u(s,q(s))= \int_s^\tau
L(\sigma, q(\sigma),\dot q(\sigma)) +\alpha d\sigma
$$
for each $s<\tau$. The existence of such a curve implies that 
$\alpha(L)=\alpha$ and that $u$ is a Weak KAM solution of $L$.
\qed

\noindent
\textbf{Proof of Proposition \ref{oneway}:}
Let us fix two points $x=(\tau,q)$ and $y=(\theta,r)$ 
in $\Tm\times M$ such that 
$h_L(x,y)=u(y)-u(x)$.
Let $q_k(t):[\tau,\theta+T_k]\lto M$ 
be a sequence of minimizing trajectories of $L$
such that $\Nm\ni T_k\lto \infty$, $q_k(\tau)=q$, 
$q_k(T_k+\theta)=r$ and 
$$\int_{\tau}^{\theta+T_k}L(t,q_k(t),\dot q_k(t))+\alpha(L) dt\lto h_L(x,y)=u(y)-u(x).
$$
Let us  consider a sequence $S_k\in [\tau,\theta+T_k]$, and assume that 
either $S_k=\tau$ (case a) or $S_k=\theta+T_k$ (case b), or both $S_k$ and $T_k-S_k$
converge to $\infty$ (case c).
By taking a subsequence, we can assume that
$S_k \mod 1$ has a limit $S$ in $\Tm$, and that the sequence 
$q_k(t+S_k)$ converge in $C^1_{loc}(I,M)$
to a limit $q(t)$ 
where $I=[0,\infty)$ in case a, $I=(-\infty,0]$ in case b, 
and $I=\Rm$ in case c.

We have, for all $s<t$ in the interior of $I$,
the equality
$$
u(t+S_k,q_k(t+S_k)))-u(s+S_k,q_k(s+S_k))=
\int_{s+S_k}^{t+S_k}
L(\sigma,q_k(\sigma),\dot q_k(\sigma))+\alpha(L) d\sigma
$$
when $k$ is large enough.
At the limit, we get 
$$
u(t+S,q(t))-u(s+S,q(s))=\int_{s+S}^{t+S}
 L(\sigma,q(\sigma-S),\dot q(\sigma-S))+
\alpha(L)d\sigma.
$$
In other words, the limit curve
 $t\lmto q(t-S)$  is calibrated by $u$ on its interval
of definition.
In case a, the limit curve
$q(t-\tau):[\tau,\infty)\lto M$
 satisfies $q(\tau-\tau)=q$. 
 This 
 implies that the point $x$ belongs to $\mI(L,u)$.

Moreover, what we have proved implies, in the terminology of the 
Appendix, that the sequence of curves 
$$
\tilde x_k(t)=(t,q_k(t),\dot q_k(t)):[\tau,\theta+T_k]\lto \Tm\times TM
$$
accumulate on $\G(L,u)$.
We conclude that $\tilde x \Cu_u \tilde y$, where
$\tilde x=\lim \tilde x(0)$,
and $\tilde y=\lim \tilde y(\theta+T_k)$.
Note that  $\tilde x=(\tau,q,\dot q(0))$, where $q(t)$
is the limit curve obtained in case a, and 
$\tilde y=(\theta,r,\dot q(0))$, where $q(t)$ is the 
limit curve obtained in case b. We have $\tilde x\in \tilde \mI(L,u)$
and $\tilde y \in \G(L,u)$.
\qed

\begin{lem}
Assume that $L_k\lto L$ and that $u_k$
is a sequence of Weak KAM solutions of $L_k$ such
that $u_k\lto u$ in $C^0$.
Then the semi-flows $\G(L_k,u_k)$ accumulate on $\G(L,u)$,
see definition in the Appendix.
If the compact sets $\G(L_k,u_k)$ converge, for the Hausdorff 
metric, to a limit $\tilde \mK$, then $\tilde \mK$ is a negatively
 invariant subset of $\G(L,u)$ such that 
 $\pi(\tilde \mK)=\Tm\times M$.
\end{lem}

\proof
It is enough to prove that, if 
$\tilde x_k(t):(-\infty,0]\lto \Tm\times TM$
is a sequence of orbits of $\G(L_k,u_k)$, then it has a subsequence
which converges uniformly to an orbit of $\G(L,u)$.
By definition of the semi-flow of $\G(L_k,u_k)$, we have
$\tilde x_k(t)=(t,q_k(t+\tau_k),\dot q_k(t+\tau_k))$, where  
$\tau_k\in \Tm$ and 
$q_k\in C^2((-\infty, \tau_k],M)$ satisfy
$$
u_k(\tau_k,q_k(\tau_k))-u_k(\tau_k-T,q_k(\tau_k-T))=\int_{\tau_k-T}^{\tau_k}
 L_k(t,q_k(t),\dot q_k(t))
+\alpha(L_k) dt
$$
for each $T\geq 0$.
We can assume that $\tau_k$ has a limit $\tau$ and that $q_k$ converges in $C^1_{loc}((-\infty,\tau[,M)$
to a curve $q(t)$.
For each $T\geq 0$ and $\epsilon>0$, we have,
for $k$ large enough,
 $$
u(\tau-\epsilon,q_k(\tau-\epsilon))-u(\tau-T,q_k(\tau-T))
=\int_{\tau-T}^{\tau-\epsilon}
 L_k(t,q_k(t),\dot q_k(t))
+\alpha(L_k) dt.
$$
By taking the limits $k\lto \infty$ and then $\epsilon\lto 0$, we get
$$
u(\tau,q(\tau))-u(\tau-T,q(\tau-T))=\int_{\tau-T}^{\tau}
 L(t,q(t),\dot q(t))
+\alpha(L) dt.
$$
This equality means that $t\lmto (t,q(t+\tau),\dot q(t+\tau))$
is an orbit of $\G(L,u)$.
If the compact sets $\G(L_k,u_k)$ converge, for the Hausdorff 
metric, to a limit $\tilde \mK$, then, what we have just done
shows that $\tilde \mK\subset \G(L,u)$.
 We have $\pi(\tilde \mK)=\Tm\times M$ because
  $\pi(\G(L_k,u_k))=\Tm\times M$.
\qed

 An interesting consequence is that, in the context above,  if 
 $\tilde x_k \in  \G(L_k,u_k)$
 and
 $\tilde y_k \in \G(L_k,u_k)$
 satisfy $\tilde x_k\lto \tilde x$, $\tilde y_k\lto \tilde y$
 and $\tilde x_k\Cu_{u_k}\tilde y_k$, then 
 $\tilde x\Cu_u \tilde y$.

\begin{lem}
Assume that $L_k\lto L$ and that $u_k$
is a sequence of Weak KAM solutions of $L_k$, such that $u_k\lto u$ in $C^0$. Assume furthermore  that 
$\tilde \mI(L_k,u_k)$ has a limit $\tilde \mK$ 
for the Hausdorff distance of 
compact sets. Then,  $\tilde \mK$ is a compact invariant subset of
 $ \tilde \mI(L,u)$. 
If $x_k\in \mI(L_k,u_k)$ and $y_k\in \mI(L_k,u_k)$ are sequences 
such that $x_k\lto x$, $y_k\lto y$ and $x_k\Cu_{u_k} y_k$,
then $x\Cu_u y$.
In addition, for each $y \in \Tm \times M$, there exists
a point $\tilde y$ above $y$ in $\G(L,u)$
and a point $\tilde x\in \tilde \mI(L,u)$
such that  $\tilde  x\Cu_u \tilde y$.
\end{lem}

\proof
We have seen that the semi-flow $\G(L_k,u_k)$ accumulates
on $\G(L,u)$.
Since $\tilde \mI(L_k,u_k)$ is invariant for this flow,
the limit $\tilde \mK$ is an invariant set of $\G(L,u)$, 
so it is contained in $\tilde \mI(L,u)$.
In order to prove the last statement, let us consider 
$y\in \Tm\times M$. There exists a point $\tilde y_k\in \G(L_k,u_k)$
above $y$.
There exists a sequence $\tilde x_k\in \tilde \mI(L_k,u_k)$
such that $\tilde x_k\Cu_{u_k} \tilde y_k$.
Taking a subsequence, we can assume that  the sequences
$\tilde y_k$ and $\tilde x_k$ have a limits
$\tilde y\in \G(L,u)$ and $\tilde x\in \tilde \mK$.
We have $\tilde x\Cu_u \tilde y$.
\qed
 \noindent
\textbf{Proof of Theorem \ref{semi}:} It  follows from:
\begin{cor}
Assume that $L_k\lto L$ and that $\tilde \mA(L_k)$ has a limit 
$\tilde \mK$
for the Hausdorff metric. Then there exists a Weak KAM solution
$u$ of $L$ such that $\tilde \mK$ is contained in the set of chain recurrent 
points of $\tilde \mI(L,u)$. If $L$ satisfies the coincidence hypothesis,
we conclude that $\tilde \mK\subset \tilde \mA(L,u)$.
\end{cor}

\proof
Let $u_k$ be a Weak KAM solution of $L_k$.
The sequence $u_k$ is equi-Lipschitz.
By adding appropriate constants, we can suppose that
it is also equi-bounded, and that  it converges to a limit $u$,
which  is  a Weak KAM solution of $L$.
Let us consider a point $\tilde x\in \tilde \mK$.
This point is the limit of a sequence 
$\tilde x_k\in \mA(L_k)\subset \tilde \mI(L_k,u_k)$.
We have $\tilde x_k \Cu_{u_k} \tilde x_k$, hence
$\tilde x\Cu_u \tilde  x$.
\qed

\section{The coincidence hypothesis}\label{coincidence}
Let us 
mention several hypotheses of non-degeneracy that appear in
the literature. The discussions in the present section are 
elaborations on the recent work of Fathi, Figalli and Rifford 
\cite{FFR}.
\begin{enumerate}
\item The quotient Aubry set is finite
\item The quotient Aubry set has 
Hausdorff dimension zero.
\item The quotient Aubry set has zero $1$-Hausdorff measure.
\item For each pair $u$, $v$ of Weak KAM solutions,
the image $(u-v)(\mA(L))$ has Lebesgue measure zero in $\Rm$.
\item For each pair $u$, $v$ of Weak KAM solutions,
the image $(u-v)(\mA(L))$ is totally disconnected in $\Rm$.
\item The coincidence hypothesis holds
\item The quotient Aubry set is totally disconnected
(which is equivalent to the statement that the static classes are the 
connected components of the Aubry set).
\end{enumerate}
The hypotheses 4. and 5. come from \cite{FFR}.
\begin{thm}
$$
1\Rightarrow 2
\Rightarrow 3 \Rightarrow 4\Leftrightarrow 5\Rightarrow 6\Rightarrow 7.
$$
\end{thm}
The implication $5\Rightarrow 4$ was pointed out to the author  by Albert Fathi.

\proof
It is obvious that 
$
1\Rightarrow 2\Rightarrow 3.
$ and that $4\Rightarrow 5$.
In order to prove that
$3 \Rightarrow 4$,
it is sufficient to notice that the difference $(u-v)$
is Lipschitz with respect to the pseudo-metric $d_L$ on $\mA(L)$.
Assuming 3, this implies that the image $(u-v)(\mA(L))$ has 
zero $1$-Hausdorff measure in $\Rm$
and therefore zero Lebesgue measure,
see \cite{FFR} for more details.
We have already seen that $6 \Rightarrow 7$.

Let us prove that
$5 \Rightarrow 6$.
The method is inspired from \cite{FFR}.
We assume 5, consider a Weak KAM solution $u$ of $L$
and two points $x$ and $y$ in $\mI(L,u)$ such that $x\Cu_u y$.
We denote by $\varphi^t$ the natural flow on $\mI(L,u)$.
We want to prove that 
$
u(y)-u(x)=h(x,y),
$
or equivalently, setting $w(z):=h(x,z)-u(z)$, that
$w(y) = w(x)$.
Contradicting this conclusion, we assume that 
$w(y)>w(x)$ (note that we always have $w(y)\geq w(x)$).
Since we assumed 5, and since $w$ is a difference of Weak KAM
solutions, there exist real numbers $a$ and $b$ such that 
$w(x)<a<b<w(y)$ and such that $w(\mA(L))\cap [a,b]=\emptyset$.
The function $w$ is non-increasing on the orbits of $\varphi^t$.
This can be seen as follows:
If $x(t)$ is an orbit of this flow, and if $s<t$, then
$$
u(x(t))-u(x(s)))= \int_s^t L(\sigma, x(\sigma), \dot x(\sigma))
+\alpha(L)d\sigma
$$
while 
$$
h_L(x,x(t))-h_L(x,x(s)))\leq
 \int_s^t L(\sigma, x(\sigma), \dot x(\sigma))
+\alpha(L)d\sigma.
$$
We conclude that the compact sets $\mI(L,u)\cap \{w\leq b\}$
and $\mI(L,u)\cap \{w\leq a\}$
are positively invariant by the flow.
It is known that each orbits of $\mI(L,u)$ is $\omega$-asymptotic
to $\mA(L)$.
All the orbits starting in $\mI(L,u)\cap \{w\leq b\}$
are thus $\omega$-asymptotic to 
$\mA(L)\cap\{w\leq b\}=\mA(L)\cap\{w< a\}$.
As a consequence, there exists $T>0$ such that 
$$
\varphi^T\big(\mI(L,u)\cap\{w\leq b\}\big)
\subset \mI(L,u)\cap\{w\leq a\}.
$$
Let us pick $\epsilon>0$ such that, for each points $z$ and $z'$ in
$\mI(L,u)$ satisfying
$w(z)\leq a$ and $d(z,z')\leq \epsilon$, we have 
$w(z')<b$.
We claim that no $(\epsilon,T)$-chain can connect $x$ and $y$.
Indeed, let  $x(t):[0,S]\lto M$ be an $(\epsilon, T)$-chain
such that $x(0)=x$.
We claim that $w(x(t))< b$ for each $t$.
Therefore, it is not possible to have $x(S)=y$.
In order to prove the claim let us denote by $\tau_i$
the jump times. We have $w(x(t))\leq w(x(\tau_i^+))$
for all $t\in [\tau_i,\tau_{i+1}[$.
So it is enough to prove that $w(x((\tau_i^+))<b$ for each $i$.
This can be proved by recurrence. If $w(x((\tau_i^+))<b$,
then, since $\tau_{i+1}\geq \tau_i+T$, we have
$$
x(\tau_{i+1}^-)\in \varphi^T\big(\mI(L,u)\cap\{w\leq b\}\big)
\subset \mI(L,u)\cap\{w \leq a\}.
$$
But then, since $d(x(\tau_{i+1}^-),\tau_{i+1}^+))\leq \epsilon$,
the way we have chosen $\epsilon$ guarantees that
 $w(x(\tau_{i+1}^+))< b$.
This ends the proof of $5 \Rightarrow 6$.

Finally, the implication  $5 \Rightarrow 4$.
follows from the next Lemma (a courtesy of Albert Fathi).
\qed

\begin{lem}\label{Fathi}
If there exist two weak KAM solutions $u$ and $v$
such that $(v-u)(\mA(L))$ has positive Lebesgue measure,
then there exists a third solution $w$ such that
$(w-u)(\mA(L))$ is a non-trivial interval.
\end{lem}

\proof
It is useful to recall that a function $w:\mA(L)\lto \Rm$
is the restriction to $\mA(L)$ of a weak KAM solution 
if and only if it satisfies
\begin{equation}\label{WK}
w(y)-w(x)\leq h_L(x,y)
\end{equation}
for each $x$ and $y$ in $\mA(L)$, see \cite{Co} and \cite{BB}.
Indeed, such a function can be extended to a weak KAM solution
on $\Tm\times M$ by the formula
$$
w(x)=\min_{a\in \mA(L)} \big( w(a)+h_L(a,x)\big).
$$
Let us denote by $A\subset \Rm$ the set $A:= (v-u)(\mA(L)),$
by $1_A(t)$ the caracteristic function of $A$, and by
$\theta_A(t)$ a primitive of $1_A(t)$.
Then, we define the function $w$ on $\mA(L)$ by
$$w:= u+\theta_A\circ (v-u).
$$
It is not hard to see that $(w-u)(\mA(L))=\theta_A(A)$
is a non-trivial interval (assuming that $A$ has positive 
Lebesgue measure).
So we have to prove that $w$ can be extended to a weak KAM
solution, or equivalentely, that (\ref{WK}) holds on $\mA(L)$.
Assume first that $(v-u)(y)\leq (v-u)(x)$.
Then, $\theta_A((v-u)(y))\leq \theta_A((v-u)(x)) $ so that
$$
w(y)-w(x)\leq u(y)-u(x)\leq h_L(x,y). 
$$
In the other case, when $(v-u)(y)\geq (v-u)(x)$,
we have, using that $\theta_A$ is $1$-Lipschitz, 
$$
w(y)-w(x)\leq u(y)-u(x)+(v-u)(y)-(v-u)(x)
=v(y)-v(x) \leq h_L(x,y).
$$
This ends the proof of Lemma \ref{Fathi}.\qed

John Mather produced in \cite{Ma:04} the example of a Lagrangian
violating 6.
However, such examples are rather exceptional.
Indeed, Ricardo Ma\~n\'e  proved 
in \cite{Mn:96} that the property of having only
one static class is generic in the following sense
(generic in the sense of Ma\~n\'e):

For each Tonelli Lagrangian $L$, there exists 
a dense $G_{\delta}$ set $\mO\subset C^{\infty}(\Tm\times M)$
of potentials
such that the property is satisfied by the Lagrangian
$L(t,x,v)+g(t,x)$ for each $g\in \mO$.

Moreover, we proved in \cite{BC} that the following property 
is generic in the sense of Ma\~n\'e (and, in a certain sense, 
satisfied outside of a singular set of infinite codimension, see
\cite{BC}):

For each $c\in H^1(M,\Rm)$, the Lagrangian $L+c$ satisfies 1.

It is even believed that the property 3 may 
hold for all smooth Lagrangians. 
The best results in that direction have
been obtained by Fathi, Figalli and Rifford in \cite{FFR},
extending earlier results of Mather, see \cite{Ma:03}
(see also \cite{So}). Their result imply that, if the dimension
of $M$ is one or two, and if $L$ is sufficiently smooth 
($C^4$ is enough), then 3 hold.
 Extending this result in higher dimension, even for analytic
  Lagrangians, is a formidable problem.

\appendix
\section{The Conley structure of flows}
This section recalls some standard facts on the topological
structure of flows on compact sets, due to Conley, see \cite{Co:88},
see also \cite{PSM}, for example,  for the extension to semi-flows.
We provide the proof of some less standard statements which are useful in the 
present paper.
It is convenient to work in an ambient metric space 
$(E,d)$.
Let $X$ be a compact subset of $E$.
A flow on $X$ is a continuous map 
$\varphi(t,x)=\varphi^t(x):\Rm\times X\lto X$ such that 
$$
\varphi^t \circ \varphi^s=\varphi^{t+s}
$$
for all $s$ and $t$ in $\Rm$.
A semi-flow on $X$ is a continuous map 
$\varphi(t,x)=\varphi^t(x):[0,\infty)\times X\lto X$ 
which 
satisfies the same relation  for $s\geq 0$ and $t\geq 0$.
We say that the  subset $Y\subset X$ is
positively  invariant by the semi-flow 
$\varphi$ if  $\varphi^t(Y)\subset Y$ for each  $t\geq 0$.
We say that $Y$ is invariant 
if it is positively invariant and if, in addition,
for each $y\in Y$ and $t\geq 0$, there exists $z\in Y$
such that $\varphi^t(z)=y$.
An $(\epsilon,T)$-chain of the semi-flow $\varphi^t$ is a piecewise continuous
curve $x(t):[0,S]\lto X$ with finitely many times of discontinuity
$S_1,\ldots, S_k\in [0,S]$ such that 
$x(t)=\varphi^{t-S_i}(x(S_i^+))$ for 
$t\in ]S_i, S_{i+1}[$, such that $S_{i+1}\geq S_i +T$ 
and such that $d(x(S_i^-),x(S_i^+))\leq \epsilon$ for each $i$.

\begin{defn}
We say that $x\Cu_{X} y$ (or, if there is no ambiguity, $x\Cu y$) if, for
 each $\epsilon>0$ and $T>0$, there exists an 
$(\epsilon,T)$-chain $x(t):[0,S]\lto X$ such that $x(0)=x$
and $x(S)=y$. 
\end{defn}

The relation $\Cu$ is closed.
It is not hard to see that we have
$$
\varphi^t(x)\Cu \varphi^s(y)
$$
for all $t\geq 0$ and $s\geq 0$ if $x\Cu y$.
This relation is satisfied for all $t$ and $s$ if $\varphi$
is a flow.

A points $x$ such that $x\Cu x$ is called chain
recurrent.
On the set of  chain recurrent points, the relation $\Cu^s$ defined by
$$
x\Cu^s y\eq \big(x\Cu y\quad\text{and}\quad y\Cu x \big) 
$$
is an equivalence relation.
Its classes of equivalence are called the chain components of $X$.
A semi-flow is said chain-recurrent if all its points are 
chain-recurrent, it is called  chain-transitive if, in addition,
 it has only
one chain component.
The chain recurrent set of a semi-flow $(X,\varphi^t)$
is contained in $X_{\infty}:= \cap_{t\geq 0} \varphi^t(X)$.
Moreover, by Lemma \ref{add} below, the chain-recurrent set of 
the semi-flow 
$(X,\varphi^t)$
is the same as the chain recurrent set of 
the restricted semi-flow $(X_{\infty},\varphi^t)$.
The following is classical (see \cite{PSM}):

\begin{prop}
Let $Y\subset X$ be the chain recurrent set of the semi-flow
 $(X,\varphi^t)$.
Then $Y$ is a compact invariant subset of $X$
which is internally chain-recurrent, which means that the
semi-flow $(Y,\varphi^t_{|Y})$ is chain recurrent.
The chain components of $(Y,\varphi^t_{|Y})$ are the chain 
components of $(X,\varphi^t)$, they are the connected components
of $Y$.
Each chain component $Z\subset Y$ is internally chain-transitive, 
in the sense that $(Z,\varphi^t_{|Z})$ is chain-transitive.
\end{prop}

\begin{defn}
Let $x_k(t):[0,T_k]\lto E$ be a sequence of (not necessarily
 continuous) maps.
 We say that $x_k$ accumulates locally uniformly on the semi-flow 
 $(X,\varphi^t)$ if 
 for each $\epsilon>0$ and $T>0$, the following property holds for infinitely many $k\in \Nm$:
 $$
  \forall S\in [0,\infty),\exists y\in X 
 \text{ s. t. }
  d(x_k(t+S),\varphi^t(y))\leq \epsilon\quad
 \forall t\in [0,T]\cap [0,T_k-S].
 $$
It is equivalent to say 
that, for each sequence $S_k\in [0,T_k]$, the sequence of curves
$x_k(t+S_k)$ has a subsequence which converges uniformly on compact 
subsets of $[0,\infty)$  to a trajectory of $\varphi^t$.

We say that a sequence $(X_k,\varphi_k^t)$ of compact sets with
 semi-flows
accumulates on $(X,\varphi^t)$ if
for each $\epsilon>0$ and $T>0$, the following property holds for infinitely many $k\in \Nm$:
 $$
 \forall x\in X_k,\exists y\in X 
 \text{ s. t. }
  d((\varphi^t_k(x),\varphi^t(y))\leq \epsilon\quad
 \forall t\in [0,T].
 $$
 It is equivalent to say that, for each $T>0$, each sequence 
 $x_k(t):[0,T]\lto X_k$ of orbits of $\varphi^t_k$ has a subsequence
 which converges uniformly to an orbit of $\varphi^t$ in  $X$.
 \end{defn}

\begin{lem}\label{compact}
Let $x_k(t):[0,T_k]\lto E$ be a sequence of 
(not necessarily continuous) curves
accumulating locally uniformly on the semi-flow $(X,\varphi^t)$.
If in addition we have $x_k(0)\lto x$ and $x_k(T_k)\lto y$,
then $x\Cu_X y$.
\end{lem}

\proof
Let us fix $\epsilon>0$ and $T>0$.
We want to prove the existence of an $(\epsilon,T)$-chain
between $x$ and $y$.
Since the curves $x_k$ accumulates  locally uniformly 
on  the flow, there exists $k\in \Nm$ such that 
$$
  \forall i\in \Nm,\exists y_i\in X 
 \text{ s. t. }
  d(x_k(t+iT),\varphi^t(y_i))\leq \epsilon/2\quad
 \forall t\in [0,T]\cap [0,T_k-iT].
 $$
 Let $y(t):[0,T_k]\lto X$ be the curve defined on 
 each interval $[iT,(i+1)T[\cap [0,T_k]$ by the expression
 $y(t)=\varphi^{t-iT}(y_i)$.
 We have 
 $$
 d(y(iT^-),y(iT^+))=
 d(\varphi^T(y_{i-1}),y_i)
 \leq d(\varphi^T(y_{i-1}),x_k(iT))+d(x_k(iT),y_i))
 \leq \epsilon
 $$
 so that the curve $y(t)$ is an appropriate pseudo-orbit.
\qed

\begin{lem}\label{conv}
Let $(X_k,\varphi^t_k)$ be a sequence  of 
compact semi-flows accumulating on $(X,\varphi^t)$.
 If $x_k\in X_k$ and $y_k\in X_k$ are sequences
 such that $x_k\lto x$, $y_k\lto y$ and 
 $x_k\Cu_{X_k}y_k$, then  $x\Cu_X y$.
\end{lem}

\proof
Let $x_k(t):[0,T_k]\lto X_k$ 
be a sequence of curves such that each $x_k(t)$
is a $(1/k,k)$-pseudo-orbit of $\varphi_k^t$ satisfying 
$x_k(0)=x_k$ and $x_k(T_k)=y_k$.
We claim that the curves $x_k$ accumulate locally uniformly on
$(X,\varphi^t)$. The conclusion then follows
from Lemma \ref{compact}.
In order to prove the claim, let us fix $\epsilon>0$ and $T>0$.
There exists $\delta\in ]0, \epsilon[$ such that, 
for all $x$ and $y$ in $X$
satisfying $d(x,y)\leq \delta$, 
and for all  $t\in [0,4T]$, we have 
$d(\varphi^t(x),\varphi^t(y))\leq \epsilon/4$.
There are infinitely many values of $k>4\max(T,1/\delta)$ such that 
$$
\forall x\in X_k,\exists y\in X 
 \text{ s. t. }
  d((\varphi^t_k(x),\varphi^t(y))\leq \delta/4\quad
 \forall t\in [0,4T].
 $$
For these values of $k$, we claim that 
$$
  \forall S\in [0,\infty),\exists y\in X 
 \text{ s. t. }
  d(x_k(t+S),\varphi^t(y))\leq \epsilon\quad
 \forall t\in [0,T]\cap [0,T_k-S].
 $$
In order to prove the claim, notice that,
given  $S\geq 0$, the curve $x_k(t+S)$ has at most one jump
on the interval $[0,T]\cap [0,T_k-S]$.
As a consequence, 
there exists a  time $\tau\in [0,T]\cap [0,T_k-S] $,
 such that 
$x_k(t+S)=\varphi^{t}_k(z)$ on $t\in [0,\tau]\cap [0,T_k-S] $,
and 
$x_k(t+S)=\varphi^{t-\tau}_k(w)$ on $t\in [\tau,T]\cap [0,T_k-S] $
for some points $z$ and $w$ in $X_k$ which satisfy 
$d(\varphi_k^{\tau}(z),w)\leq 1/k\leq \delta/4$.
There exist two points $Z$ and $W$ in $ X$ such that
 $$
 d((\varphi^t_k(z),\varphi^t(Z))\leq \delta/4\quad
 \forall t\in [0,2T]
 $$
and 
$$d((\varphi^t_k(w),\varphi^t(W))\leq \delta/4\quad
 \forall t\in [0,2T].
 $$
In particular, we have 
$$d(\varphi^{\tau}(Z),W)
\leq d(\varphi^{\tau}(Z),\varphi^{\tau}_k(z))+
d(\varphi^{\tau}_k(z),w)+d(w,W)\leq \delta.
$$
Let us set $y:=\varphi^{-\tau}(Z)$.
For each $t\in [0,\tau]$, we have
$$
d(x_k(t+S),\varphi^t(Z))=d(\varphi^t_k(z),\varphi^t(Z))
\leq \delta/4\leq \epsilon,
$$
and for each $t\in [\tau,T]\cap [0,T_k-S] $, we have
$$
d(x_k(t+S),\varphi^t(Z))=d(\varphi^{t-\tau}_k(w),
\varphi^{t-\tau}(\varphi^{\tau}(Z)))
$$
$$
\leq d(\varphi^{t-\tau}(w),\varphi^{t-\tau}(W))+
d(\varphi^{t-\tau}(W),\varphi^{t-\tau}(\varphi^{\tau}(Z)))
\leq \delta/4+\epsilon/4\leq \epsilon.
$$
This proves the claim.
\qed

The following obvious remark will be useful:
\begin{lem}
Assume that the semi-flows $(X_k,\varphi^t_k)$ accumulate on
 $(X,\varphi)$,
and that $X_k$ converge to $Y$ for the Hausdorff metric.
Then $Y$ is a compact positively invariant subset of $X$, and 
$(X_k,\varphi^t_k)$ accumulate on $(Y,\varphi^t)$.
\end{lem}

If $(X,\varphi^t)$ is a semi-flow, then,
for $k\in \Nm$, the set $X_k:=\varphi^k(X)$
is positively invariant, so that $(X_k,\varphi^t)$
is itself a semi-flow. It is not hard to see that
the relations $\Cu_X$ and $\Cu_{X_k}$ coincide on $X_k$.
Defining now 
$
X_{\infty}:=\cap_{k\in \Nm} X_k
$
we have:

\begin{lem}\label{add}
The set $X_{\infty}$ is positively invariant;
the relations $\Cu_{X}$ and $\Cu_{X_{\infty}}$ coincide on 
$X_{\infty}$.
\end{lem}
\proof
The relations $\Cu_{X_k}$, $k\in \Nm$ all coincide on $X_{\infty}$.
We have to prove that
$$x\Cu_{X_{k}} y \Rightarrow x\Cu_{X_{\infty}} y.
$$
This follows from Lemma \ref{conv} since 
the semi-flows $(X_k,\varphi^t)$ accumulate on 
$(X_{\infty}, \varphi^t)$.
\qed

\small
\bibliographystyle{amsplain}

\begin{thebibliography}{30}

\bibitem{ETDS} P. Bernard:
\textit{Homoclinic orbits to invariant sets of quasi-integrable exact maps.} Ergodic Theory Dynam. Systems
\textbf{ 20} (2000), no. 6, 1583--1601. 

\bibitem{Fourier}  P. Bernard
\textit{  Connecting orbits of time dependent 
Lagrangian systems.}
 Ann. Inst. Fourier {\bf  52} (2002), 1533--1568.


\bibitem{ENS}
P. Bernard, 
\textit{Existence of $C^{1,1}$ critical sub-solutions of the Hamilton-Jacobi equation on compact manifolds,}
Ann. Sci. E. N. S. (4) \textbf{40} (2007), No. 3, 445-452 

\bibitem{AMS}  P. Bernard,
\textit{The dynamics of pseudographs in convex Hamiltonian 
systems,} Journ. AMS, \textbf{21} (2008) 615--669.


\bibitem{BB}
P. Bernard, B. Buffoni:
\textit{Weak KAM pairs and Monge-Kantorovich duality}
Advanced Studies in Pure Mathematics, \textbf{47}, No. 2, (2007) 397-420,
 Asymptotic analysis and singularity.

\bibitem{BC}
P. Bernard, G. Contreras:
\textit{Generic properties of families of Lagrangian systems,} 
Annals of Maths \textbf{167} (2008) No 3.

\bibitem{Co:88}
C. Conley:\textit{ The gradient structure of a flow. I.}
 With a comment by R. Moeckel.  
 Ergodic Theory Dynam. Systems  \textbf{8}  (1988),  Charles Conley Memorial Issue, 11--26,




\bibitem{Co}
G. Contreras:
\textit{ Action potential and weak KAM solutions,}
 Calc. Var. Partial Differential Equations 
 \textbf{13} (2001), no. 4, 427--458.



\bibitem{Mn:96}
R. Ma\~n\'e,  
\textit{Generic properties and problems of minimizing measures of
   {L}agrangian systems}, Nonlinearity \textbf{9} (1996), no.~2, 273--310.


\bibitem{Mn:97}
R. Ma\~n\'e,  
\textit{ Lagrangian flows: The dynamics of
globally minimizing orbits,} Bol. Soc. Bras. Mat, \textbf{28}
(1997) 141-153


\bibitem{Ma:91}
J. N. Mather, 
\textit{ Action minimizing invariant
measures for positive definite Lagrangian systems,}
 Math.\
Z. \textbf{207} 169--207 (1991)


\bibitem{Ma:93}
J. N. Mather, 
\textit{ Variational construction of connecting orbits, }
Ann. Inst. Fourier, \textbf{43} (1993), 1349-1368.



\bibitem{Ma:03}
J. N. Mather,
\textit{Total disconnectedness of the quotien Aubry set
in low dimensions}
Comm. Pure App. Math. \textbf{56} (2003), 1178--1183.


\bibitem{Ma:04}
J. N. Mather, \textit{Examples of Aubry sets}, Ergod. Th. \& Dyn. Syst.
 \textbf{24} (2004) 1667--1723.

\bibitem{FFR}
A. Fathi, A. Figalli and L. Rifford,
\textit{On the Hausdorff dimension of the Mather quotient.}

\bibitem{Fa:book}
A. Fathi, Weak KAM Theorem in Lagrangian Dynamics,
 Book to appear.



\bibitem{PSM}
M. Patr\~ao, L. A. B  San Martin,
\textit{Semiflows on Topological Spaces: Chain Transitivity and semigroups,}
Jour. Dyn. Diff. Eqn. \textbf{19} (2007) no. 1 155--180.


\bibitem{So}
A. Sorrentino,
\textit{On the total disconnectedness of the quotient Aubry set}
Ergodic Theory Dynam. Systems. \textbf{28} (2008) no. 1 267--290.

\end{thebibliography}
\providecommand{\bysame}{\leavevmode\hbox
to3em{\hrulefill}\thinspace}

\end{document}